# THE TWO-TYPE RICHARDSON MODEL WITH UNBOUNDED INITIAL CONFIGURATIONS

By Maria Deijfen and Olle Häggström

*Stockholm University and Chalmers University of Technology*

The two-type Richardson model describes the growth of two competing infections on $\mathbb{Z}^d$ and the main question is whether both infection types can simultaneously grow to occupy infinite parts of $\mathbb{Z}^d$. For bounded initial configurations, this has been thoroughly studied. In this paper, an unbounded initial configuration consisting of points $x = (x_1, \ldots, x_d)$ in the hyperplane $\mathcal{H} = \{x \in \mathbb{Z}^d : x_1 = 0\}$ is considered. It is shown that, starting from a configuration where all points in $\mathcal{H} \setminus \{\mathbf{0}\}$ are type 1 infected and the origin $\mathbf{0}$ is type 2 infected, there is a positive probability for the type 2 infection to grow unboundedly if and only if it has a strictly larger intensity than the type 1 infection. If, instead, the initial type 1 infection is restricted to the negative $x_1$-axis, it is shown that the type 2 infection at the origin can also grow unboundedly when the infection types have the same intensity.

**1. Introduction.** One of the simplest models for spatial growth and competition is the Richardson model, introduced in Richardson [13]. The original version describes the growth of a single infectious entity on $\mathbb{Z}^d$, but the mechanism can be extended to comprise two entities, making it a model for competition on $\mathbb{Z}^d$; see Häggström and Pemantle [7, 8]. This paper is concerned with the two-type version of the model in $d \geq 2$ dimensions, started from a configuration where one of the entities occupies one single site in an infinite "sea" of the other entity.

The dynamics of the one-type Richardson model are such that an uninfected site becomes infected at a rate proportional to the number of infected nearest neighbors and, once infected, it never recovers. This is equivalent to first-passage percolation with i.i.d. exponential passage times. The main result, dating back to Richardson [13] and Kesten [10], is that, for bounded initial configurations, the infection grows linearly in time in each direction and, on the scale $1/t$, the set of infected points at time $t$ converges almost surely to a deterministic shape; see Theorem 2.1.

In the two-type version of the model, a second infection type is introduced and the two infections, referred to as type 1 and type 2, respectively, grow simultaneously on $\mathbb{Z}^d$, the dynamics being that an uninfected site becomes type $i$ infected at a rate proportional to the number of type $i$ infected nearest neighbors and then stays (type $i$) infected forever. The model has two parameters, denoted by $\lambda_1$ and $\lambda_2$, indicating the intensities of the infections.

For the two-type Richardson model, the main challenge lies in deciding if it is possible for both infection types to simultaneously grow to occupy infinite parts of $\mathbb{Z}^d$. We will denote by $G_i$ the event that the type $i$ infection reaches sites arbitrarily far away from the origin. It has been conjectured that in the two-type Richardson model in $d \geq 2$ dimensions, started from one single site of each infection type, the event $G_1 \cap G_2$ has positive probability if and only if $\lambda_1 = \lambda_2$; see Häggström and Pemantle [7]. The if direction of the conjecture was proved for $d = 2$ by Häggström and Pemantle [7] and, in the general case, independently by Garet and Marchand [6] and Hoffman [9].

As for the "only if" direction of the conjecture, the best result to date is that, with one of the intensities fixed, mutual unbounded growth has probability 0 for all but at most countably many values of the other intensity; see Häggström and Pemantle [8]. The full result would indeed follow from this if it could be proven that the probability of the event $G_1 \cap G_2$ is monotone, in the sense that it decreases as the difference between the intensities increases. This monotonicity is certainly believed to hold on the $\mathbb{Z}^d$-lattice. However, it has turned out to be very hard to prove and, in fact, there are other graphs where it actually fails; see Deijfen and Häggström [2]. On the other hand, in Deijfen and Häggström [3], it is shown that, as long as we restrict ourselves to bounded initial configurations on $\mathbb{Z}^d$, the particular choice of initial sets is not important in deciding whether the event of mutual unbounded growth for the two infection types has positive probability or not.

The purpose of the present paper is to study the two-type Richardson model with *unbounded* initial configurations. The question is as before: can both types simultaneously infect infinitely many sites? With infinite initial configurations for both types, the answer is (other than in absurd cases) obviously "yes." We therefore restrict to cases where type 1 starts with infinitely many sites and type 2 with finitely many (for concreteness, only one) and the question of coexistence then reduces to that of whether type 2 can survive. Write $(x_1, \ldots, x_d)$ for the coordinates of a point $x \in \mathbb{Z}^d$ and define $\mathcal{H} = \{x : x_1 = 0\}$ and $\mathcal{L} = \{x : x_1 \leq 0 \text{ and } x_i = 0 \text{ for all } i \geq 2\}$. Writing

**0** for the origin, the initial configurations that we will consider are

(1)
$$\begin{aligned} I(\mathcal{H}): \quad &\text{all points in } \mathcal{H}\backslash\{\mathbf{0}\} \text{ are type 1 infected and} \\ &\mathbf{0} \text{ is type 2 infected;} \\ I(\mathcal{L}): \quad &\text{all points in } \mathcal{L}\backslash\{\mathbf{0}\} \text{ are type 1 infected and} \\ &\mathbf{0} \text{ is type 2 infected.} \end{aligned}$$

Interestingly, it turns out that the set of parameter values that allows for coexistence is slightly different for these two configurations. Let $P_{\mathcal{H},\mathbf{0}}^{\lambda_1,\lambda_2}$ (resp. $P_{\mathcal{L},\mathbf{0}}^{\lambda_1,\lambda_2}$) denote the probability measure associated with a two-type process started from configuration $I(\mathcal{H})$ (resp. $I(\mathcal{L})$) and note that, by time scaling, we may restrict our attention to the case $\lambda_1 = 1$. Our main result is as follows.

THEOREM 1.1. *For the Richardson model in $d \geq 2$ dimensions, we have:*

(a) $P_{\mathcal{H},\mathbf{0}}^{1,\lambda}(G_2) > 0$ *if and only if* $\lambda > 1$;
(b) $P_{\mathcal{L},\mathbf{0}}^{1,\lambda}(G_2) > 0$ *if and only if* $\lambda \geq 1$.

In words, when the type 2 infection is strictly stronger than the type 1 infection, there is a positive probability for the type 2 infection to grow unboundedly in both configurations. Intuitively, the type 2 infection can use its higher intensity to rush away along the $x_1$-axis and achieve an unsurpassable lead over the type 1 infection. When the infection types have the same intensity, the type 2 infection can grow unboundedly from the configuration $I(\mathcal{L})$, but not from $I(\mathcal{H})$, where the initial disadvantage for the type 2 infection turns out to be too severe.

Before proceeding, we mention the following question, which has been pointed out to us by Itai Benjamini as well as an anonymous referee: Consider the case when the infections have the same intensity (i.e., $\lambda = 1$) and assume that the type 1 infection occupies not only the negative $x_1$-axis (as in $I(\mathcal{L})$), but a cone around it. What is the critical slope of the cone for which there is a positive probability for the type 2 infection at the origin to grow unboundedly? It seems likely, as suggested to us by Itai Benjamini, that the critical case is when the cone fills the whole left half-space. That the type 2 infection cannot survive when the type 1 infection occupies the whole left half-space follows from Theorem 1.1(a). That infinite type 2 growth has positive probability for any smaller type 1 cone remains to be proven and goes beyond the scope of this paper.

The rest of the paper is organized as follows. In Section 2, some results are described for a one-type process started with the entire hyperplane $\mathcal{H}$ infected at time 0. These results are then used in Sections 3, 4 and 5 to prove the claims in Theorem 1.1 for $\lambda > 1$, $\lambda = 1$ and $\lambda < 1$ respectively; these cases will be referred to as *supercritical*, *critical* and *subcritical*. Finally, some concluding remarks appear in Section 6.

**2. The one-type process started from $\mathcal{H}$.** In this section, we recall from the literature the results needed to deduce that the asymptotic speed along the $x_1$-axis of the growth in a one-type process started with all sites in $\mathcal{H}$ infected at time 0 is the same as the asymptotic speed along the axes in a one-type process started from the origin. We also state a result on the rate of convergence of the speed, originating with Kesten [11], and introduce a hampered version of the process where only edges in a certain part of $\mathbb{Z}^d$ can be used.

First, consider a one-type process with intensity $\lambda$ started from the origin. In such a process, each edge $e$ of the $\mathbb{Z}^d$-lattice has an exponential random variable $\tau(e)$ with parameter $\lambda$ associated with it. The travel time of a path $\Gamma$ is defined as

$$T(\Gamma) := \sum_{e \in \Gamma} \tau(e)$$

and, for $x \in \mathbb{Z}^d$, the time when the point $x$ is infected is given by

$$T^{\mathbf{0}}(x) := \inf\{T(\Gamma) : \Gamma \text{ is a path from } \mathbf{0} \text{ to } x\}.$$

Write $\mathbf{n} = (n, 0, \ldots, 0)$. A basic result then, obtained from the subadditive ergodic theorem, is the existence of a constant $\mu_\lambda > 0$ such that $T^{\mathbf{0}}(\mathbf{n})/n \to \mu_\lambda$ almost surely and in $L_1$ as $n \to \infty$; see, for example, Kingman [12]. Defining $\mu = \mu_1$, it follows from a simple time-scaling argument that $\mu_\lambda = \mu \lambda^{-1}$, hence we have

$$(2) \qquad \lim_{n \to \infty} \frac{T^{\mathbf{0}}(\mathbf{n})}{n} = \mu \lambda^{-1} \qquad \text{a.s. and in } L_1.$$

More generally, for any $x \in \mathbb{Z}^d$ with $x \neq \mathbf{0}$, there is a constant $\mu(x) > 0$ such that $T^{\mathbf{0}}(nx)/n \to \mu(x)\lambda^{-1}$ as $n \to \infty$. Hence, the infection grows linearly in each fixed direction and the asymptotic speed of growth is an almost sure constant. That this is also true when all directions are considered simultaneously is stated in the shape theorem, which is the main result for the model. To formulate it, write $\xi^{\mathbf{0}}(t)$ for the set of infected points at time $t$ in a process started from the origin and let $\bar{\xi}^{\mathbf{0}}(t) \subset \mathbb{R}^d$ be a continuum version of $\xi^{\mathbf{0}}(t)$ obtained by replacing each $x \in \xi^{\mathbf{0}}(t)$ by a unit cube centered at $x$.

THEOREM 2.1. *There exists a compact convex set $A$ such that, for any $\varepsilon > 0$, almost surely*

$$(1-\varepsilon)\lambda A \subset \frac{\bar{\xi}^{\mathbf{0}}(t)}{t} \subset (1+\varepsilon)\lambda A$$

*for large $t$.*

An "in probability" version of this result was established by Richardson [13] and strengthened to the present form by Kesten [13]. Note that $A$ is the unit ball in the norm defined by $\mu(x)$, that is, $A = \{x : \mu(x) \leq 1\}$.

Now, consider the Richardson model started with all sites in the hyperplane $\mathcal{H}$ infected. We will show that the asymptotic speed of growth in the direction of the $x_1$-axis in such a process is, in fact, the same as the speed of growth along the axes in a process started with only the origin infected. To this end, write $T^{\mathcal{H}}(x)$ for the time when the point $x \in \mathbb{Z}^d$ is infected in a process started from $\mathcal{H}$.

PROPOSITION 2.2. *In the unit rate one-type process, we have, as $n \to \infty$, that $T^{\mathcal{H}}(\mathbf{n})/n \to \mu$ in $L_1$.*

PROOF. Let $\mathcal{H}_n$ be the hyperplane at $x_1$-coordinate $n$, that is, $\mathcal{H}_n = \{x \in \mathbb{Z}^d : x_1 = n\}$. (In this notation, we have $\mathcal{H} = \mathcal{H}_0$.) The first important observation is that, in the first-passage percolation representation of the model that is used in this paper, the time when the point $\mathbf{n}$ is infected in a process started from $\mathcal{H}$ is the same as the time when the first (in time) point belonging to $\mathcal{H}$ is infected in a process started from $\mathbf{n}$. Furthermore, the minimal travel time from $\mathbf{n}$ to $\mathcal{H}$ has the same distribution as the minimal travel time from the origin to $\mathcal{H}_n$, which we denote by $T^{\mathbf{0}}(\mathcal{H}_n)$. Hence, we have that

$$(3) \qquad \frac{T^{\mathcal{H}}(\mathbf{n})}{n} \stackrel{d}{=} \frac{T^{\mathbf{0}}(\mathcal{H}_n)}{n},$$

where $\stackrel{d}{=}$ denotes equality in distribution. Just as in establishing (2), the subadditive ergodic theorem can be applied to show that $T^{\mathbf{0}}(\mathcal{H}_n)/n$ converges in $L_1$ (and almost surely) to a constant $c$. Using (3), the proposition follows if we can show that $c = \mu$. Clearly, since $\mathbf{n} \in \mathcal{H}_n$, we have $T^{\mathbf{0}}(\mathcal{H}_n) \leq T^{\mathbf{0}}(\mathbf{n})$, which gives $c \leq \lim_n T^{\mathbf{0}}(\mathbf{n})/n = \mu$. The reverse inequality follows from the fact that the asymptotic shape of a process started from the origin is convex. Indeed, having $c < \mu$ would contradict the convexity of the asymptotic shape $A$ stipulated in Theorem 2.1. Hence, $c = \mu$, as desired. $\square$

Next, we state a result on the convergence rate in (2). The result originates with Kesten [11] and is formulated here for the time $T^{\mathbf{0}}(\mathcal{H}_n)$ when the hyperplane $\mathcal{H}_n$ is reached by the infection. In the original formulation, the estimate concerns the time $T^{\mathbf{0}}(\mathbf{n})$ when the single point $\mathbf{n}$ is infected, but it is pointed out that the bound also applies to passage times to hyperplanes. Kesten [11] also contains results on the convergence rate in the shape theorem; related results, improving some of the bounds of Kesten, can be found in Alexander [1].

THEOREM 2.3 (Kesten [11]). *There exist constants $c_1$, $c_2$ and $c_3$ such that*

$$\mathbf{P}\left(\left|\frac{T^{\mathbf{0}}(\mathcal{H}_n) - \mathbf{E}[T^{\mathbf{0}}(\mathcal{H}_n)]}{\sqrt{n}}\right| \geq x\right) \leq c_1 e^{-c_2 x} \qquad \text{for } x \leq c_3 n.$$

Combining this estimate with Proposition 2.2 gives the following lemma, which will be useful in controlling the type 1 infection in a two-type process started according to $I(\mathcal{H})$.

LEMMA 2.4. *In a unit rate one-type process, for any $\varepsilon > 0$, there exist constants $c$ and $c'$ such that*

(4) $$\mathbf{P}(T^{\mathcal{H}}(\mathbf{n}) \leq (1-\varepsilon)n\mu) \leq ce^{-c'\sqrt{n}}$$

*for large $n$.*

PROOF. Trivially, we have

$$\mathbf{P}(T^{\mathcal{H}}(\mathbf{n}) \leq (1-\varepsilon)n\mu) \leq \mathbf{P}(|T^{\mathcal{H}}(\mathbf{n}) - n\mu| \geq \varepsilon n\mu).$$

To obtain $|T^{\mathcal{H}}(\mathbf{n}) - n\mu| \geq \varepsilon n\mu$, by the triangle inequality, at least one of the quantities $|T^{\mathcal{H}}(\mathbf{n}) - \mathbf{E}[T^{\mathcal{H}}(\mathbf{n})]|$ and $|\mathbf{E}[T^{\mathcal{H}}(\mathbf{n})] - n\mu|$ must exceed $\varepsilon n\mu/2$. By Proposition 2.2, we will, indeed, not have $|\mathbf{E}[T^{\mathcal{H}}(\mathbf{n})] - n\mu| \geq \varepsilon n\mu/2$ when $n$ is large, hence

$$\mathbf{P}(|T^{\mathcal{H}}(\mathbf{n}) - n\mu| \geq \varepsilon n\mu) \leq \mathbf{P}\left(|T^{\mathcal{H}}(\mathbf{n}) - \mathbf{E}[T^{\mathcal{H}}(\mathbf{n})]| \geq \frac{\varepsilon n\mu}{2}\right)$$

for large $n$. As observed in the proof of Proposition 2.2, the passage time $T^{\mathcal{H}}(\mathbf{n})$ has the same distribution as $T^{\mathbf{0}}(\mathcal{H}_n)$. Choosing $x = \sqrt{n}\varepsilon\mu/2$ (which is clearly smaller than $c_3 n$ when $n$ is large) in Theorem 2.3, we hence obtain

$$\mathbf{P}\left(|T^{\mathcal{H}}(\mathbf{n}) - \mathbf{E}[T^{\mathcal{H}}(\mathbf{n})]| \geq \frac{\varepsilon n\mu}{2}\right) \leq c_1 e^{-\sqrt{n}\varepsilon\mu c_2/2}$$

and the lemma follows. □

Using Lemma 2.4, it turns out that we can establish that the convergence in Proposition 2.2 also holds in the almost sure sense.

PROPOSITION 2.5. *As $n \to \infty$, we have that $T^{\mathcal{H}}(\mathbf{n})/n \to \mu$ almost surely.*

PROOF. We need to show that for any $\varepsilon > 0$,

(5) $$\limsup_{n \to \infty} \frac{T^{\mathcal{H}}(\mathbf{n})}{n} \leq (1+\varepsilon)\mu$$

almost surely and that

(6) $$\liminf_{n\to\infty} \frac{T^{\mathcal{H}}(\mathbf{n})}{n} \geq (1-\varepsilon)\mu$$

almost surely. First, note that (2) implies that (5) holds with $T^{\mathbf{0}}(\mathbf{n})$ in place of $T^{\mathcal{H}}(\mathbf{n})$. But, obviously, $T^{\mathcal{H}}(\mathbf{n}) \leq T^{\mathbf{0}}(\mathbf{n})$, so (5) is established.

Next, in order to show (6), fix $N < \infty$ in such a way that (4) holds for all $n \geq N$. From Lemma 2.4, we then have that the expected number of $n \geq 0$ for which the event $\{T^{\mathcal{H}}(\mathbf{n}) \leq (1-\varepsilon)n\mu\}$ happens is at most

$$N + \sum_{n=N}^{\infty} ce^{-c'\sqrt{n}},$$

which is finite. Hence, by the Borel–Cantelli lemma, (6) is established. □

The last consideration before moving on to the two-type process is to show that the growth of a one-type process restricted to a "tube" around the $x_1$-axis behaves approximately as an unrestricted process when the tube is large. This will be needed to control the type 2 infection at the origin in a two-type process started from the configuration $I(\mathcal{H})$. To formulate the result, we introduce a new, hampered version of the one-type process by placing "walls" in $\mathbb{Z}^d$ restricting the growth in all directions except one. More precisely, we consider a process with the same dynamics as the original one, but where only sites in the set

(7) $$\Omega_b := \{x \in \mathbb{Z}^d : |x_i| \leq b \text{ for all } i \neq 1\}$$

are susceptible to the infection. We write $\xi^{b*}(t)$ for the set of infected points at time $t$ in such a process, started with a single infection at the origin $\mathbf{0}$. The following lemma says that $\Omega_b$ is filled with infection linearly in time.

LEMMA 2.6. *Consider a hampered one-type process with rate $\lambda$. For any dimension $d$, there is a real number $\mu_{\lambda,b} > 0$ such that, for any $\varepsilon \in (0, \mu_{\lambda,b}^{-1})$, almost surely*

$$\{x \in \Omega_b : |x_1| \leq (1-\varepsilon)t\mu_{\lambda,b}^{-1}\} \subset \xi^{b*}(t) \subset \{x \in \Omega_b : |x_1| \leq (1+\varepsilon)t\mu_{\lambda,b}^{-1}\}$$

*for all sufficiently large $t$.*

Being a completely standard adaptation of the proof of the shape theorem, the proof of this lemma is omitted.

The constant $\mu_{\lambda,b}$ is the analog of $\mu\lambda^{-1}$ in the unhampered process. That is, if $T^{b*}(x)$ denotes the time when the point $x$ is infected in a hampered process, then we have $\mu_{\lambda,b} = \lim_n T^{b*}(\mathbf{n})/n$. When $b$ is large, it is reasonable to expect that the speed of growth in the hampered process is close to the speed of an unhampered process and hence that $\mu_{\lambda,b}$ is close to $\mu\lambda^{-1}$ for large $b$. This intuition is confirmed by the next lemma.

LEMMA 2.7. *As $b \to \infty$, we have that $\mu_{\lambda,b} \to \mu\lambda^{-1}$.*

The proof of this lemma is again a straightforward adaptation, this time of the proof of Lemma 4.4 in Deijfen, Häggström and Bagley [4], where the same result is established for a continuum counterpart of the Richardson model. It is therefore omitted.

**3. The supercritical case.** Our main task in this section is to prove coexistence when $\lambda > 1$.

PROPOSITION 3.1. *For any $\lambda > 1$ and any $d \geq 2$, we have*
$$P^{1,\lambda}_{\mathcal{H},\mathbf{0}}(G_2) > 0.$$

Here and later, we will make use of the following convenient way of constructing the two-type Richardson model: attach to each edge $e$ in the edge set $E_{\mathbb{Z}^d}$ of the $\mathbb{Z}^d$-lattice independently two independent exponential variables $\tau_1(e)$ and $\tau_2(e)$ with mean 1 (resp. $\lambda^{-1}$), indicating the time it takes for the type 1 (resp. type 2) infection to traverse the edge.

PROOF OF PROPOSITION 3.1. Fix $\lambda > 1$. For a hampered one-type process in a region $\Omega_b$, as in (7), Lemma 2.7 tells us that $\lim_{b\to\infty} \mu_{\lambda,b} \to \mu\lambda^{-1}$ and, since $\mu\lambda^{-1} < \mu$, we have $\mu_{\lambda,b} < \mu$ when $b$ is sufficiently large. Fix such a $b$ and define
$$\Omega_b^+ = \{x \in \Omega_b : x_1 \geq 1\}.$$

For $x \in \mathbb{Z}^d$, define $S^{1,\mathcal{H}\setminus\{\mathbf{0}\}}(x)$ as $\inf_\Gamma \sum_{e \in \Gamma} \tau_1(e)$, where the infimum is over all paths starting at $\mathcal{H} \setminus \{\mathbf{0}\}$ and ending at $x$. Furthermore, for $x \in \Omega_b^+$, define $S_b^{2,\mathbf{0}}(x)$ as $\inf_\Gamma \sum_{e\in\Gamma} \tau_2(e)$, where, this time, the infimum is over all paths starting at $\mathbf{0}$ *and passing through vertices in $\Omega_b$ only*.

The point of these definitions is the following observation, which is easy to see and which turns out to be instrumental in proving Proposition 3.1.

> Suppose that at least one vertex in $\Omega_b^+$ is eventually infected by type 1 and let $x$ be the first (in time) vertex in $\Omega_b^+$ for which this happens. We then have
> $$S^{1,\mathcal{H}\setminus\{\mathbf{0}\}}(x) \leq S_b^{2,\mathbf{0}}(x).$$

So, if we can show that

(8) $\qquad \mathbf{P}(S^{1,\mathcal{H}\setminus\{\mathbf{0}\}}(x) > S_b^{2,\mathbf{0}}(x) \text{ for all } x \in \Omega_b^+) > 0,$

then we know that, with positive probability, all $x \in \Omega_b^+$ eventually become infected by type 2, and the proposition follows.

To do this, define $\mu^* = \frac{\mu+\mu_{\lambda,b}}{2}$ so that $\mu_{\lambda,b} < \mu^* < \mu$ and note that, by Proposition 2.5, we have the existence of some (random) $M_1 < \infty$ such that

$$S^{1,\mathcal{H}\setminus\{\mathbf{0}\}}(x) > \mu^* x_1 \qquad \text{for all } x \in \Omega_b^+ \text{ with } x_1 > M_1.$$

Likewise, Lemma 2.6 guarantees the existence of some (again random) $M_2 < \infty$ such that

$$S_b^{2,\mathbf{0}}(x) < \mu^* x_1 \qquad \text{for all } x \in \Omega_b^+ \text{ with } x_1 > M_2.$$

Taking $M = \max\{M_1, M_2\}$ yields

$$S^{1,\mathcal{H}\setminus\{\mathbf{0}\}}(x) > S_b^{2,\mathbf{0}}(x) \qquad \text{for all } x \in \Omega_b^+ \text{ with } x_1 > M,$$

which is very close to proving (8). To rigorously get from here to (8), we employ the following conditioning argument. Choose an $m$ such that

(9) $\qquad \mathbf{P}(S^{1,\mathcal{H}\setminus\{\mathbf{0}\}}(x) > S_b^{2,\mathbf{0}}(x) \text{ for all } x \in \Omega_b^+ \text{ with } x_1 > m) > 0$

and define $E_b^{+,m}$ as the set of edges $\langle x, y \rangle \in E_{\mathbb{Z}^d}$ such that either

$$x = \mathbf{0} \quad \text{and} \quad y = \mathbf{1} = (1, 0, \ldots, 0)$$

or

$$x, y \in \Omega_b^+, x_1 \leq m, y_1 \leq m+1.$$

We write $D_m$ for the event in (9). We will condition on the $\tau_1(e)$ variables for all edges $e \in E_{\mathbb{Z}^d}$, together with the $\tau_2(e)$ variables for all edges $e \in E_{\mathbb{Z}^d} \setminus E_b^{+,m}$. Let $A$ be the event that

(10) $\qquad \mathbf{P}(D_m | \{\tau_1(e)\}_{e \in E_{\mathbb{Z}^d}}, \{\tau_2(e)\}_{e \in E_{\mathbb{Z}^d} \setminus E_b^{+,m}}) > 0$

and note that since $\mathbf{P}(D_m) > 0$, we must also have $\mathbf{P}(A) > 0$. On the event $A$, there exists a random $\gamma > 0$ (depending on $\{\tau_1(e)\}_{e \in E_{\mathbb{Z}^d}}$ and $\{\tau_2(e)\}_{e \in E_{\mathbb{Z}^d} \setminus E_b^{+,m}}$) such that if $\tau_2(e) < \gamma$ for all $e \in E_b^{+,m}$, then $D_m$ occurs. By further decreasing $\gamma > 0$, we can also ensure (due to the fact that $S^{1,\mathcal{H}\setminus\{\mathbf{0}\}}(x) > 0$ for all $x$ almost surely) that if $\tau_2(e) < \gamma$ for all $e \in E_b^{+,m}$, then we also have

$$S^{1,\mathcal{H}\setminus\{\mathbf{0}\}}(x) > S_b^{2,\mathbf{0}}(x) \qquad \text{for all } x \in \Omega_b^+ \text{ with } x_1 \leq m.$$

Finally, note that on the event $A$—that is, under the conditioning in (10)—the probability that $\tau_2(e) < \gamma$ for all $e \in E_b^{+,m}$ is strictly positive; this is simply because the edge set $E_b^{+,m}$ is finite and $\gamma > 0$. Hence,

$$\mathbf{P}(S^{1,\mathcal{H}\setminus\{\mathbf{0}\}}(x) > S_b^{2,\mathbf{0}}(x) \text{ for all } x \in \Omega_b^+ | A) > 0.$$

Since $\mathbf{P}(A) > 0$, (8) follows and the proof is complete. $\square$

The other result concerning the supercritical case that we need to prove is the following.

COROLLARY 3.2. *For any $\lambda > 1$ and any $d \geq 2$, we have*

$$P^{1,\lambda}_{\mathcal{L},\mathbf{0}}(G_2) > 0.$$

One way to prove this, is to note that the proof of Proposition 3.1 can be easily adapted to handle the corollary. Or, to be a bit more careful, we can invoke the proposition itself, together with the following easy lemma. For disjoint subsets $\xi_1$ and $\xi_2$ of $\mathbb{Z}^d$, write $P^{\lambda_1,\lambda_2}_{\xi_1,\xi_2}$ for the law of the two-type process with parameters $\lambda_1$ and $\lambda_2$, with all sites in $\xi_1$ initially infected by type 1, all in $\xi_2$ by type 2 and all others uninfected.

LEMMA 3.3. *Suppose that $\xi_1$ and $\xi_2$ are disjoint subsets of $\mathbb{Z}^d$, and $\xi'_1$ and $\xi'_2$ likewise. If $\xi_1 \subseteq \xi'_1$ and $\xi_2 \supseteq \xi'_2$, then*

$$P^{\lambda_1,\lambda_2}_{\xi_1,\xi_2}(G_2) \geq P^{\lambda_1,\lambda_2}_{\xi'_1,\xi'_2}(G_2).$$

PROOF. Couple the two processes using the same $\tau_1(e)$ and $\tau_2(e)$ variables. Writing $(\xi_1(t), \xi_2(t))$ for the state of the first process at time $t$ in the obvious way and similarly for the second process, it is straightforward to show that the relations $\xi_1(t) \subseteq \xi'_1(t)$ and $\xi_2(t) \supseteq \xi'_2(t)$ are preserved as $t$ increases. Letting $t \to \infty$ proves the lemma. □

PROOF OF COROLLARY 3.2. Define $\mathcal{L}'$ in the same way as $\mathcal{L}$, except with the roles of the $x_1$- and $x_2$-coordinates interchanged, that is, let

$$\mathcal{L}' = \{x : x_2 \leq 0 \text{ and } x_i = 0 \text{ for all } i \neq 2\}.$$

Then, by symmetry, we have $P^{1,\lambda}_{\mathcal{L}',\mathbf{0}}(G_2) = P^{1,\lambda}_{\mathcal{L},\mathbf{0}}(G_2)$. But, since $\mathcal{L}' \subset \mathcal{H}$, we can invoke Lemma 3.3 in order to deduce that $P^{1,\lambda}_{\mathcal{L}',\mathbf{0}}(G_2) \geq P^{1,\lambda}_{\mathcal{H},\mathbf{0}}(G_2)$, the latter probability being positive by Proposition 3.1. □

**4. The critical case.** The critical case $\lambda = 1$ is, in some ways, the most interesting, particularly in view of the fact that whether or not a single type-2 infection at the origin has a chance to survive against an infinite army of type-1 enemies depends on whether these initially get to occupy all of $\mathcal{H} \setminus \{\mathbf{0}\}$ or just $\mathcal{L} \setminus \{\mathbf{0}\}$. The two results we need to prove in this section are the following.

PROPOSITION 4.1. *For any $d \geq 2$, we have $P^{1,1}_{\mathcal{H},\mathbf{0}}(G_2) = 0$.*

PROPOSITION 4.2. *For any $d \geq 2$, we have $P^{1,1}_{\mathcal{L},\mathbf{0}}(G_2) > 0$.*

For $\lambda = 1$, it is often convenient to replace the construction following Proposition 3.1 by one that involves not two but just one travel time variable $\tau(e)$ associated with each $e \in E_{\mathbb{Z}^d}$. Here, the $\tau(e)$'s are taken to be i.i.d. exponentials with mean 1 and represent the time it takes for either of the two infection types to traverse the edge. We will employ this construction in the proofs of both Propositions 4.1 and 4.2. A nice feature of the construction is the following. Suppose that the process starts at time 0 with the nodes in $\xi_1 \subset \mathbb{Z}^d$ infected by type 1 and the nodes in $\xi_2 \subset \mathbb{Z}^d$ infected by type 2. We write $T^{\xi_1}(x)$ for $\inf_\Gamma \sum_{e \in \Gamma} \tau(e)$, where the infimum is over all paths starting in $\xi_1$ and ending at $x$, and define $T^{\xi_2}(x)$ analogously. Then, $x$ becomes infected precisely at time $\min\{T^{\xi_1}(x), T^{\xi_2}(x)\}$ and, furthermore,

(11) $\qquad x$ gets infected by type 1 if and only if $T^{\xi_1}(x) < T^{\xi_2}(x)$.

The following lemma (an easy variation of Lemma 3.3) will be useful in the proof of Proposition 4.1.

LEMMA 4.3. *Consider the (one-type) Richardson model with parameter 1 on a bounded degree graph $L$ with vertex set $V$ and edge set $E$ starting at time 0 with the set $\xi \subset V$ initially infected. Also, consider the same model with parameter 1 on another graph $L' = (V', E')$ starting at time 0 with the set $\xi' \subset V'$ initially infected. If $L$ is a subgraph of $L'$, in the sense that $V \subseteq V'$ and $E \subseteq E'$, and, furthermore, we have $\xi \subseteq \xi'$, then, for any $t > 0$ and any $\eta \subset E$, the probability that all vertices in $\eta$ are infected at time $t$ is no greater for the process on $L$ than for the process on $L'$.*

PROOF. Simply couple the two processes in such a way that for each $e \in E$, they use the same $\tau(e)$ variable. $\square$

PROOF OF PROPOSITION 4.1. We will show that

(12) $\qquad \mathbf{P}(\text{type 2 infects infinitely many sites in the half-space } \{x : x_1 \geq 1\}) = 0.$

Once that is done, we can, by symmetry, infer the corresponding statement for the other half-space and thus conclude that $\mathbf{P}(G_2) = 0$.

Recall, for integer $b$, that $\mathcal{H}_b$ is the set of vertices $x \in \mathbb{Z}^d$ whose $x_1$-coordinate is $b$ (so that, in particular, $\mathcal{H}_0 = \mathcal{H}$). For $x \in \mathbb{Z}^d$ and $y \in \mathcal{H}$, we write $y \to x$ for the event that the infimum $T^{\mathcal{H}}(x)$ is obtained by a path from $y$ to $x$. Beginning with all nodes in $\mathcal{H}$ infected, we think of $y \to x$ as meaning that the infection eventually hitting $x$ descends from $y$. Since the $\tau(e)$'s are independent with a continuous distribution, we have, for fixed $x$, that the event $y \to x$ happens for only one $y$ almost surely (i.e., there are no ties). Thus,

(13) $$\sum_{y \in \mathcal{H}} \mathbf{P}(y \to x) = 1.$$

Fix $b \geq 1$, $y, y' \in \mathcal{H}$ and $x, x' \in \mathcal{H}_b$ in such a way that $x$ and $y$ differ only in their first coordinate, likewise for $x'$ and $y'$. Symmetry implies that

$$\mathbf{P}(y \to x') = \mathbf{P}(y' \to x)$$

and, in conjunction with (13), this implies that

$$\sum_{x \in \mathcal{H}_b} \mathbf{P}(y \to x) = 1 \tag{14}$$

for any $y \in \mathcal{H}$. Let $X_b$ denote the number of sites $x \in \mathcal{H}_b$ such that $\mathbf{0} \to x$. Then, to prove (12) is the same as showing that

$$\mathbf{P}\left(\sum_{b=1}^{\infty} X_b = \infty\right) = 0. \tag{15}$$

By (14), we have $\mathbf{E}[X_b] = 1$, whence $\mathbf{P}(X_b = \infty) = 0$ for any $b$, so in order to prove (15), it is enough to show that

$$\mathbf{P}\left(\lim_{b \to \infty} X_b = 0\right) = 1. \tag{16}$$

Ideally, since $\mathbf{E}[X_b] = 1$ for each $b$, we would now like to endow the sequence $(X_1, X_2, \ldots)$ with a martingale structure (with respect to some filtration). Since a nonnegative martingale converges almost surely and this one presumably could not converge to anything other than 0, that would settle (16). But we cannot see how to do this (it is probably not even possible) and will settle for a different solution that, although a bit less clean, is still reminiscent of a martingale approach.

To this end, we first define, for $x \in \mathbb{Z}^d$ with $x_1 \leq b$ and $y \in \mathcal{H}$, the passage time $T_b^y(x)$ as $\inf_\Gamma \sum_{e \in \Gamma} \tau(e)$, where the infimum is over all paths from $y$ to $x$ that *do not pass through any vertex $z$ with $z_1 > b$*. Also, define

$$T_b^{\mathcal{H}}(x) = \inf_{y \in \mathcal{H}} T_b^y(x) \tag{17}$$

and for $y' \in \mathcal{H}$, let $y' \xrightarrow{b} x$ denote the event that the infimum in (17) is attained for $y = y'$. Finally, write $X_b^*$ for the number of vertices in $\mathcal{H}_b$ such that $\mathbf{0} \xrightarrow{b} x$. By the same argument as for $X_b$, we have, for any $b \geq 1$, that

$$\mathbf{E}[X_b^*] = 1. \tag{18}$$

We now claim, crucially, that

$$\{X_b \geq 1\} \text{ implies } \{X_b^* \geq 1\}. \tag{19}$$

To see this, assume that $X_b \geq 1$ and take $x$ to be the vertex among those in $\mathcal{H}_b$ satisfying $\mathbf{0} \to x$ for which $T^{\mathbf{0}}(x)$ is smallest. Then, the path from $\mathbf{0}$ to $x$ cannot pass through any other vertex $z$ in $\mathcal{H}_b$ (because if that were the case, $z$ would satisfy $\mathbf{0} \to z$ with a smaller $T_z^{\mathbf{0}}$). Hence, $T_b^{\mathbf{0}}(x) = T^{\mathbf{0}}(x)$. Since $T_b^y(x) \geq T^y(x)$ for all $y \in \mathcal{H}$, it follows that $\mathbf{0} \xrightarrow{b} x$ and the claim (19) is warranted.

Now, write $G_2^*$ for the event that $\limsup_{b \to \infty} X_b^* > 0$. Using (19), we can show that (16)—and thereby the proposition—follows if we can show that

$$\mathbf{P}(G_2^*) = 0, \tag{20}$$

so this is what we set out to prove.

Our argument will involve the filtration $\{\mathcal{F}_b\}_{b=0}^{\infty}$, where $\mathcal{F}_b$ is the $\sigma$-field generated by the $\tau(e)$-variables for all edges $e = \langle x, y \rangle$ with $x_1, y_1 \leq b$. Note that, for $x \in \mathcal{H}_b$, the event $\mathbf{0} \xrightarrow{b} x$ and the random variables $T_b^{\mathcal{H}}(x)$, $T_b^{\mathbf{0}}(x)$ and $X_b^*$ are all $\mathcal{F}_b$-measurable. Lévy's 0-1 law (see, e.g., Durrett [5]) tells us that

$$\lim_{b \to \infty} \mathbf{P}(G_2^* | \mathcal{F}_b) = \mathbb{I}_{G_2^*} \tag{21}$$

almost surely, where $\mathbb{I}_{G_2^*}$ is the indicator of the event $G_2^*$. The conclusion of the argument will be to show that, with probability 1, $\mathbf{P}(G_2^* | \mathcal{F}_b)$ does not converge to 1.

Fix $\varepsilon > 0$ small, in such a way that $n = \frac{2}{\varepsilon}$ is an integer. For $x, y \in \mathbb{Z}^d$, write $\text{dist}(x, y)$ for the graph-theoretic distance between $x$ and $y$ in the $\mathbb{Z}^d$ lattice. We claim that there exists a $u < \infty$ such that, for any $b$, the event $D_{b,u}$ has probability at most $\frac{\varepsilon}{2}$, where we define $D_{b,u}$ as the event that

$X_b^* \leq n$, while, for some $x, y \in \mathcal{H}_b$ such that $\mathbf{0} \xrightarrow{b} x$ and $\text{dist}(x, y) \leq n$, we have $T_b^{\mathcal{H}}(y) \geq T_b^{\mathcal{H}}(x) + u$.

To see this, we will invoke a comparison with the classical one-type Richardson model on $\mathbb{Z}^{d-1}$ with infection rate 1 and starting at time 0 with a single infection at the origin. Since, almost surely, any finite set of sites is eventually infected in this model, we can find a $u$ such that the probability that all sites within distance $n$ from the origin are infected by time $u$ is at least $1 - \frac{\varepsilon}{2n}$.

Returning to the model in $\mathbb{Z}^d$, consider (for some fixed $k \in \{1, \ldots, n\}$) the vertex $x$ that has the $k$th smallest value of $T_b^{\mathcal{H}}(x)$ among those for which $\mathbf{0} \xrightarrow{b} x$ (provided at least $k$ such vertices exist). For a fixed realization of the process up to time $T_b^{\mathcal{H}}(x)$, a stochastic domination argument now shows that the probability that all sites in $\mathcal{H}_b$ within distance $n$ from $x$ are infected by time $T_b^{\mathcal{H}}(x) + u$ is at least $1 - \frac{\varepsilon}{2n}$; to see this, apply Lemma 4.3 with $L$ equal to $\mathcal{H}_b$ with edges between Euclidean nearest neighbors (this graph is isomorphic to $\mathbb{Z}^{d-1}$) and $\xi = \{x\}$, $L'$ equal to $\mathbb{Z}^d$ restricted to vertices $z$ with $z_1 \leq b$ and $\xi'$ equal to the set of vertices infected by time $T_b^{\mathcal{H}}(x)$ in the original process.

Summing the resulting complementary probability bound $\frac{\varepsilon}{2n}$ from 1 to $n$ gives the desired bound

$$\mathbf{P}(D_{b,u}) \leq \frac{\varepsilon}{2}.$$

Furthermore, Markov's inequality applied to (18) yields $\mathbf{P}(X_b^* > n) \leq \frac{1}{n} = \frac{\varepsilon}{2}$. If we now define the event

$$\widetilde{D}_{b,u} = D_{b,u} \cup \{X_b^* > n\},$$

we obtain

$$\mathbf{P}(\widetilde{D}_{b,u}) \leq \frac{\varepsilon}{2} + \frac{\varepsilon}{2} = \varepsilon.$$

Let $\neg$ denote set complement. The next crucial claim is that (given $\varepsilon > 0$ as above) there exists a $\delta > 0$ (independent of $b$) such that the event $\neg \widetilde{D}_{b,u}$ implies

(22) $$\mathbf{P}(X_{b+1}^* = 0 | \mathcal{F}_b) \geq \delta.$$

To see this, assume that $\neg \widetilde{D}_{b,u}$ occurs and consider the following event $A_{b,u}$ which, together with $\neg \widetilde{D}_{b,u}$, is enough to guarantee that $X_{b+1}^* = 0$. Namely, let $A_{b,u}$ be the event that:

- for all $e = \langle x, y \rangle$ with $y \in \mathcal{H}_{b+1}$, $x \in \mathcal{H}_b$, $\mathbf{0} \not\xrightarrow{b} x$ and $x$ within distance $n$ from some $z \in \mathcal{H}_b$ such that $\mathbf{0} \xrightarrow{b} z$, we have $\tau(e) \leq 1$;
- for all $e = \langle x, y \rangle$ with $x, y \in \mathcal{H}_{b+1}$ within distance $n+1$ from some $z \in \mathcal{H}_b$ such that $\mathbf{0} \xrightarrow{b} z$, we have $\tau(e) \leq 1$;
- for all $e = \langle x, y \rangle$ with $y \in \mathcal{H}_{b+1}$, $x \in \mathcal{H}_b$ and $\mathbf{0} \xrightarrow{b} x$, we have $\tau(e) \geq u+n+1$.

The point of this definition is that if $\neg \widetilde{D}_{b,u}$ and $A_{b,u}$ occur, then none of the vertices $z \in \mathcal{H}_b$ with $\mathbf{0} \xrightarrow{b} z$ will have time to infect their neighbor in $\mathcal{H}_{b+1}$ before infection creeps in from another direction, and therefore $X_{b+1}^*$ will equal 0. These requirements on $A_{b,u}$ are guide restrictive, but note that, on the event $\neg \widetilde{D}_{b,u}$, there are at most $n$ edges $e$ for which we require that $\tau(e) \geq u+n+1$ and there are at most $nd(2n)^{d-1}$ edges $e$ for which we require that $\tau(e) \leq 1$. Thus, on the event $\neg \widetilde{D}_{b,u}$, we have

$$\mathbf{P}(A_{b,u}|\mathcal{F}_b) \geq (e^{-(u+n+1)})^n (1-e^{-1})^{nd(2n)^{d-1}},$$

which is a small number indeed, but strictly positive, and the whole point of this exercise is that it does not depend on $b$. Thus, we have for any $b$ on the event $\neg \widetilde{D}_{b,u}$ that (22) holds with $\delta = (e^{-(u+n+1)})^n (1-e^{-1})^{nd(2n)^{d-1}}$. Since $X_{b+1}^* = 0$ precludes $G_2^*$, we also have on the event $\neg \widetilde{D}_{b,u}$ that

$$\mathbf{P}(G_2^*|\mathcal{F}_b) \leq 1 - \delta.$$

So, in order for the limit $\lim_{b\to\infty} \mathbf{P}(G_2^*|\mathcal{F}_b)$ in (21) to equal 1, we need $\widetilde{D}_{b,u}$ to occur for all sufficiently large $b$. But $\tilde{D}_{b,u}$ was defined in such a way as to guarantee that $\mathbf{P}(\widetilde{D}_{b,u}) \leq \varepsilon$ for any $b$, so we can conclude that $\mathbf{P}(G_2^*) \leq \varepsilon$. Since $\varepsilon > 0$ could be taken to be arbitrarily small, we obtain $\mathbf{P}(G_2^*) = 0$ and the proof is complete. $\square$

That was rather involved. Fortunately, the proof of Proposition 4.2 is somewhat more straightforward.

PROOF OF PROPOSITION 4.2. Recalling the notation $\mathbf{n} = (n, 0, \ldots, 0)$, the initial configuration in Proposition 4.2 consists of a single type-2 infection at $\mathbf{0}$, competing against type-1 infections at $-\mathbf{1}, -\mathbf{2}, \ldots$. Intuitively, the best hope for type 2 seems to be to rush off along the positive $x_1$-axis. Assuming this to be the case and viewing the model in terms of the $\tau(e)$-variables as before, we will set out to prove that

(23) $$\liminf_{n\to\infty} \mathbf{P}\left(T^{\mathbf{0}}(\mathbf{n}) = \inf_{\mathbf{m}\leq 0} T^{\mathbf{m}}(\mathbf{n})\right) > 0.$$

Note that the event in (23) is precisely the event that the vertex $\mathbf{n}$ becomes infected by type 2. Hence, if we can prove (23), we can deduce that, with positive probability, infinitely many vertices on the positive $x_1$-axis are infected and the proposition will follow.

Symmetry implies that for any $n \geq 1$, we have

(24) $$\mathbf{P}\left(T^{\mathbf{0}}(\mathbf{n}) = \inf_{\mathbf{m}\leq 0} T^{\mathbf{m}}(\mathbf{n})\right) = \mathbf{P}\left(T^{\mathbf{0}}(\mathbf{n}) = \inf_{l\geq n} T^{\mathbf{0}}(\mathbf{l})\right).$$

We will work with the right-hand side of (24), the advantage of this being that it has a useful interpretation in terms of the one-type Richardson model. Namely, for the one-type model starting at time 0 with a single infected site at the origin, the event $\{T^{\mathbf{0}}(\mathbf{n}) = \inf_{l\geq n} T^{\mathbf{0}}(\mathbf{l})\}$ is precisely the event that the node $\mathbf{n}$ is infected before any node further away on the positive $x_1$-axis is infected.

To this end, for the one-type Richardson model, define $Y(t)$ as the number of nodes on the positive $x_1$-axis that are infected by time $t$. Further, define $Y^{\rightarrow}(t)$ as the number of nodes on the positive $x_1$-axis that are infected by time $t$, with the additional property that at the time of their infection, they became the rightmost infected node on the $x_1$-axis.

We have, from (2), that

$$\lim_{t\to\infty} \frac{Y(t)}{t} = \mu^{-1}$$

almost surely, where $\mu > 0$ is the time constant discussed in Section 2. What about $\frac{Y^{\rightarrow}(t)}{t}$? At any time $t$, there exists a rightmost infected node on the $x_1$-axis and this node infects its neighbor-to-the-right at rate 1. Every time such an infection occurs, $Y^{\rightarrow}(t)$ increases by 1. Hence, the process $\{Y^{\rightarrow}(t)\}_{t\geq 0}$ stochastically dominates a rate-1 Poisson process, whence

$$\liminf_{t\to\infty} \frac{Y^{\rightarrow}(t)}{t} \geq 1 \tag{25}$$

almost surely. Again by (2), we have, with probability 1, that for any $\varepsilon > 0$, eventually all infected nodes on the positive $x_1$-axis at time $t$ have an $x_1$-coordinate that does not exceed $(1+\varepsilon)t\mu^{-1}$, that is, the number of infected nodes on the positive $x_1$-axis at time $t$ does not exceed $(1+\varepsilon)t\mu^{-1}$. In conjunction with (25), this implies that

$$\liminf_{n\to\infty} n^{-1} \sum_{j=1}^{n} \mathbb{I}_{\{T^{\mathbf{0}}(\mathbf{j})=\inf_{l\geq j} T^{\mathbf{0}}(\mathbf{l})\}} \geq \mu$$

almost surely. Hence, by Fatou's lemma, we have

$$\liminf_{n\to\infty} n^{-1} \sum_{j=1}^{n} \mathbf{P}\Big(T^{\mathbf{0}}(\mathbf{j}) = \inf_{l\geq j} T^{\mathbf{0}}(\mathbf{l})\Big) \geq \mu,$$

implying that

$$\liminf_{n\to\infty} \mathbf{P}\Big(T^{\mathbf{0}}(\mathbf{n}) = \inf_{l\geq n} T^{\mathbf{0}}(\mathbf{l})\Big) \geq \mu.$$

Using the identity (24), this implies (23), and the proof is complete. □

**5. The subcritical case.** For the subcritical case $\lambda < 1$, we shall see that the following result holds, which is more general than the subcritical cases of both part (a) and part (b) of Theorem 1.1.

PROPOSITION 5.1. *Consider the two-type Richardson model on $\mathbb{Z}^d$, $d \geq 2$, with types 1 and 2 having respective intensities 1 and $\lambda$, starting with type 1 in $\xi_1 \subset \mathbb{Z}^d$ and type 2 in $\xi_2 \subset \mathbb{Z}^d$. If $\xi_1$ is infinite, $\xi_2$ is finite and $\lambda < 1$, then the event $G_2$ of unbounded survival for type 2 has probability 0.*

It turns out that this result is a direct consequence of the following proposition which was instrumental in proving the main result in Häggström and Pemantle [8] mentioned in Section 1. Let $A$ be the asymptotic shape for the Richardson model, as defined by Theorem 2.1.

PROPOSITION 5.2. *For any $\lambda < 1$ and any $\varepsilon > 0$, we have*

$$\lim_{r\to\infty} \sup_{\xi_1,\xi_2} P^{1,\lambda}_{\xi_1,\xi_2}(G_2) = 0,$$

*where the supremum is over all initial configurations $(\xi_1,\xi_2)$ such that*

(26)
$\xi_2$ *is contained in $rA$, while*

$\xi_1$ *is not contained in $(1+\varepsilon)rA$.*

In fact, the proposition as stated in Häggström and Pemantle [8] dealt only with the case where $\xi_1$ was finite, but the generalization to infinite $\xi_1$ follows immediately from Lemma 3.3. Proposition 5.1 now follows from Proposition 5.2 upon noting that if $\xi_1$ is infinite and $\xi_2$ is finite, then the pair $(\xi_1, \xi_2)$ satisfies (26) for all sufficiently large $r$.

**6. Concluding remarks.** To see that our main result, Theorem 1.1, has now been proven is just a matter of collecting the results from the previous three sections. Theorem 1.1(a) follows from Propositions 3.1 (supercritical case), 4.1 (critical case) and 5.1 (subcritical case), while Theorem 1.1(b) follows from Corollary 3.2 (supercritical case) and Propositions 4.2 (critical case) and 5.1 (subcritical case).

We end the paper with the observation that Proposition 4.2 allows us to construct a simple proof of the fact that infinite coexistence is possible in the critical case $\lambda = 1$ starting from finitely many infected nodes of each type (recall the result of Deijfen and Häggström [3] that the particular choice of finite initial configurations does not matter, as long as type 1 does not already "strangle" type 2 or vice versa). There already exist several proofs of this result—Häggström and Pemantle [7] produced one for $d = 2$, while Garet and Marchand [6] and Hoffman [9] proved the result for arbitrary $d$—but, since the result is central to the study of the two-type Richardson model, we feel it is worth the extra effort to state a new, simple proof.

THEOREM 6.1. *For the critical ($\lambda = 1$) two-type Richardson model on $\mathbb{Z}^d$ in any dimension $d \geq 2$, there exists an $n$ such that if the model starts with a single type-1 infection at $\mathbf{0}$ and a single type-2 infection at $\mathbf{n}$, then infinite coexistence has positive probability.*

PROOF. Consider the usual edge representation of the critical model, where each edge is assigned a $\tau(e)$ representing the time it takes either infection to traverse it. Also, as usual, for $\xi \subset \mathbb{Z}^d$ and $x \in \mathbb{Z}^d$, write $T^\xi(x)$ for the sum of the $\tau(e)$'s along the fastest path starting in $\xi$ and ending at $x$. Define two random sequences $\{X_n\}_{n=-\infty}^\infty$ and $\{Y_n\}_{n=-\infty}^\infty$ as follows. Set

$$X_n = \begin{cases} 1, & \text{if } T^{\mathbf{n}}(z) < T^{\{\ldots,\mathbf{n-3},\mathbf{n-2},\mathbf{n-1}\}}(z) \text{ for infinitely many } z \in \mathbb{Z}^d, \\ 0, & \text{otherwise,} \end{cases}$$

and

$$Y_n = \begin{cases} 1, & \text{if } T^{\mathbf{n}}(z) < T^{\{\mathbf{n+1},\mathbf{n+2},\mathbf{n+3},\ldots\}}(z) \text{ for infinitely many } z \in \mathbb{Z}^d, \\ 0, & \text{otherwise.} \end{cases}$$

Let $a = \mathbf{P}(X_0 = 1)$. Proposition 4.2 tells us that $a > 0$. The process $\{X_n\}_{n=-\infty}^\infty$ is stationary, so $\mathbf{P}(X_n = 1) = a$ for any $n$. By symmetry, $\mathbf{P}(Y_n = 1) = a$ also holds for any $n$. Furthermore, $\{X_n\}_{n=-\infty}^\infty$ arises in a stationary way from an i.i.d. process and is therefore ergodic, so $\mathbf{P}(X_n = 1 \text{ for some } n \geq 1) = 1$. Hence, we can find an $n$ such that

$$\mathbf{P}(Y_0 = 1, X_n = 1) > 0. \tag{27}$$

On the event $\{Y_0 = 1, X_n = 1\}$, we have (by definition of the two processes) that $T^{\mathbf{0}}(z) < T^{\mathbf{n}}(z)$ for infinitely many $z \in \mathbb{Z}^d$ and that $T^{\mathbf{n}}(z) < T^{\mathbf{0}}(z)$ for infinitely many $z \in \mathbb{Z}^d$. Thus, (27) guarantees that infinite coexistence has positive probability for the two-type model starting at $\mathbf{0}$ and at $\mathbf{n}$. $\square$

## REFERENCES

[1] ALEXANDER, K. (1993). A note on some rates of convergence in first-passage percolation. *Ann. Appl. Probab.* **3** 81–90. MR1202516
[2] DEIJFEN, M. and HÄGGSTRÖM, O. (2006a). Nonmonotonic coexistence regions for the two-type Richardson model on graphs. *Electron. J. Probab.* **11** 331–344. MR2217820
[3] DEIJFEN, M. and HÄGGSTRÖM, O. (2006b). The initial configuration is irrelevant for the possibility of mutual unbounded growth in the two-type Richardson model. *Combin. Probab. Comput.* **15** 345–353. MR2216472
[4] DEIJFEN, M., HÄGGSTRÖM, O. and BAGLEY, J. (2004). A stochastic model for competing growth on $\mathbb{R}^d$. *Markov Process. Related Fields* **10** 217–248. MR2082573
[5] DURRETT, R. T. (1991). *Probability: Theory and Examples*. Wadsworth and Brooks/Cole, Pacific Grove, CA. MR1068527
[6] GARET, O. and MARCHAND, R. (2005). Coexistence in two-type first-passage percolation models. *Ann. Appl. Probab.* **15** 298–330. MR2115045
[7] HÄGGSTRÖM, O. and PEMANTLE, R. (1998). First passage percolation and a model for competing spatial growth. *J. Appl. Probab.* **35** 683–692. MR1659548
[8] HÄGGSTRÖM, O. and PEMANTLE, R. (2000). Absence of mutual unbounded growth for almost all parameter values in the two-type Richardson model. *Stochastic Process. Appl.* **90** 207–222. MR1794536
[9] HOFFMAN, C. (2005). Coexistence for Richardson type competing spatial growth models. *Ann. Appl. Probab.* **15** 739–747. MR2114988
[10] KESTEN, H. (1973). Discussion contribution. *Ann. Probab.* **1** 903. MR0356192
[11] KESTEN, H. (1993). On the speed of convergence in first-passage percolation. *Ann. Appl. Probab.* **3** 296–338. MR1221154
[12] KINGMAN, J. F. C. (1968). The ergodic theory of subadditive stochastic processes. *J. Roy. Statist. Soc. Ser. B* **30** 499–500. MR0254907
[13] RICHARDSON, D. (1973). Random growth in a tessellation. *Proc. Cambridge Phil. Soc.* **74** 515–528. MR0329079

DEPARTMENT OF MATHEMATICS  
STOCKHOLM UNIVERSITY  
106 91 STOCKHOLM  
SWEDEN  
E-MAIL: mia@math.su.se

DEPARTMENT OF MATHEMATICS  
CHALMERS UNIVERSITY OF TECHNOLOGY  
412 96 GOTHENBURG  
SWEDEN  
E-MAIL: olleh@math.chalmers.se